\documentclass[11pt,reqno]{amsart}
\usepackage[top=1.25in,bottom=1.25in,left=1in,right=1in]{geometry}
\usepackage{amsfonts, amsmath, amssymb, amscd, amsthm, graphicx}
\usepackage{enumitem}
\usepackage[final, 
colorlinks=true]{hyperref}
\usepackage{verbatim}

\usepackage{tikz,tikz-cd,amscd}
\usetikzlibrary{calc}
\usepackage{tikz,mathdots,color}
\usetikzlibrary{decorations.pathreplacing,decorations.markings}
\usetikzlibrary{shapes,shapes.arrows,positioning}

\tikzset{->-/.style={decoration={
  markings,
  mark=at position #1 with {\arrow{>}}},postaction={decorate}}}
\tikzset{->-/.default=.5}

\theoremstyle{plain}
\newtheorem{theorem}{Theorem}
\newtheorem{lemma}[theorem]{Lemma}
\newtheorem{proposition}[theorem]{Proposition}

\newtheorem{claim}[theorem]{Claim}
\newtheorem*{claim*}{Claim}

\theoremstyle{definition}
\newtheorem{definition}[theorem]{Definition}
\newtheorem{notation}[theorem]{Notation}

\newtheorem{conjecture}[theorem]{Conjecture}

\numberwithin{equation}{section}
\numberwithin{theorem}{section}

\newcommand{\RR}{\mathbb{R}}

\newcommand{\calB}{\mathcal{B}}
\newcommand{\calC}{\mathcal{C}}

\newcommand{\calM}{\mathcal{M}}

\newcommand{\calR}{\mathcal{R}}

\newcommand{\bu}{\mathbf{u}}
\newcommand{\bv}{\mathbf{v}}

\newcommand{\One}{\mathbf{1}}

\newcommand{\hell}{\hat{\ell}}

\newcommand{\oF}{\overline{F}}

\newcommand{\entropy}{\mathfrak{h}}  
\newcommand{\entropyinf}{\entropy_{\rm inf}}  
\newcommand{\entropysup}{\overline{\entropy}}  
    
\newcommand{\abs}[1]{\left\lvert {#1} \right\rvert}

\newcommand{\from}{\colon\thinspace}

\newcommand{\param}%
	{{\mathchoice{\mkern1mu\mbox{\raise2.2pt\hbox{$\centerdot$}}\mkern1mu}%
	{\mkern1mu\mbox{\raise2.2pt\hbox{$\centerdot$}}\mkern1mu}%
	{\mkern1.5mu\centerdot\mkern1.5mu}{\mkern1.5mu\centerdot\mkern1.5mu}}}

\DeclareMathOperator{\rk}{rk}
\DeclareMathOperator{\spec}{spec}
\DeclareMathOperator{\Out}{Out}

\begin{document}

\title{Subgraph Entropy}


\author{Tarik Aougab}
\address{Department of Mathematics, Haverford College, Haverford, PA 19041, USA}
\curraddr{}
\email{taougab@haverford.edu}

\author{Matt Clay}
\address{Department of Mathematical Sciences, University of Arkansas, Fayeteville, AR 72701, USA}
\curraddr{}
\email{mattclay@uark.edu}

\author{Tawfiq Hamed}
\address{Birzeit University, West Bank, Palestine}
\curraddr{}
\email{tawfiq9009@gmail.com}

\begin{abstract} 

Given $r \geq 3$, we prove that there exists $\lambda >0$ depending only on $r$ so that if $G$ is a metric graph of rank $r$ with metric entropy $1$, then there exists a proper subgraph $H$ of $G$ with metric entropy at least $\lambda$. This answers a question of the second two authors together with Rieck. We interpret this as a graph theoretic version of the Bers Lemma from hyperbolic geometry, and explain some connections to the pressure metric on the Culler-Vogtmann Outer Space. 

\end{abstract}

\date{\today}

\maketitle

\section{Introduction}

For any $g \geq 2$, there is a constant $\beta_g$ such that for any closed hyperbolic surface $X$ there is a collection of $3g-3$ curves on $X$ that are pairwise non-isotopic where each has length at most $\beta_g$. The existence of such a constant $\beta_g$ was first shown by Bers \cite{ar:Bers74}, and the minimal such number $\beta_g$ is called the \emph{Bers constant}.  This constant plays a key role in understanding the topology of moduli space, the geometry of the mapping class group and of hyperbolic $3$-manifolds, and the structure of Teichm{\"u}ller space (\cite{Farb2012}, \cite{ar:M96},  \cite{ar:MM00}, \cite{ar:Brock03}). Other important estimates include the work of Buser \cite{th:Buser80}, Buser-Sepp{\"a}l{\"a} \cite{ar:BS92}, and Schmutz \cite{ar:schmutz94} on the growth of the Bers constant and related functions, such as the maximal length of a systole as a function of the topology of a given hyperbolic surface.  The best known bound on the Bers constant is due to Parlier~\cite{ar:Parlier24}. 

Over the last 40 years, techniques from the study of hyperbolic surfaces have played an influential role in studying the moduli spaces of metric graphs and the outer automorphism group of a free group (\cite{ar:Bestvina}, \cite{ar:BestvinaFeighn}, \cite{ar:CullerVogtmann}, \cite{ar:Kao}, \cite{ar:PS14}, \cite{ar:ACR23}). For example, the second two authors and Rieck \cite{ar:ACR23} studied a pair of piecewise Riemannian metrics on the Culler-Vogtmann Outer Space that were inspired by the \textit{thermodynamical interpretation} of the Weil-Petersson metric on Teichm{\"u}ller space, due to Wolpert, McMullen, and Thurston (\cite{ar:Wol}, \cite{ar:Mc}). It was Thurston who first defined a Riemannian metric on Teichm{\"u}ller space by considering a limit of Hessians associated to a sequence of ``randomly'' chosen closed geodesics; later, Wolpert confirmed that this was a multiple of the Weil-Petersson metric \cite{ar:Wol}, and McMullen helped to develop a more general framework for the so-called \textit{thermodynamic formalism} for making sense of this phenomenon \cite{ar:Mc}. It is therefore no surprise that the metrics studied in \cite{ar:ACR23} are intimately related to counting closed cycles on objects being parameterized by the Outer space, namely, (marked) metric graphs.

In the spirit of porting over intuition from hyperbolic geometry to the setting of graphs and free groups, the main goal of this paper is to derive an analogue of the Bers constant for metric graphs. Of course, depending on how one defines the normalized Outer Space, the most naive translations of Bers' theorem to the setting of graphs is either trivial, or can not possibly hold. We can see this even in the simplest case of a rose $R_n$ on $n$ petals; indeed, perhaps one is interested in finding a maximal collection of pairwise disjoint (except at the basepoint) cycles on $R_n$ each with length at most some $\beta_n$. If the normalized Outer Space is defined to consist of metric graphs of total volume $1$, this is trivially true (and one can choose $\beta_n = 1$) and not particularly useful. On the other hand, the metrics on Outer Space mentioned in the previous paragraph require that one chooses a different normalization-- one that is suited to considering the growth of the collection of all cycles as a function of length-- namely, selecting metric graphs with unit \textit{entropy}, the exponential growth rate of the number of reduced cycles. And since there are roses with unit entropy and arbitrarily large volume, no version of Bers' theorem can hold. 

Note that we have already stumbled upon a key difference between the geometry of hyperbolic surfaces and metric graphs: any closed hyperbolic surface \textit{automatically} has unit entropy and has area depending only on $g$. Since one way of deriving Bers' theorem is to use the Gauss-Bonnet theorem and the bound it produces on total area, one might hope there is a way of reframing the argument in terms of entropy instead of area. Indeed, one can prove the following, \textit{and} that it is equivalent to the existence of the Bers constant $\beta_{g}$ (see the Section \ref{sec:prelim} for a proof and for a more formal statement, and for all relevant definitions):

\begin{theorem} \label{thm:BersForGraphs}  Given $\chi$ and $L$, there is a constant $\eta >0$ such that if $X$ is a hyperbolic surface with Euler characteristic $\chi$ and with $\ell(\partial X) \le L$, there is a proper essential subsurface $Y \subsetneq X$ with entropy at least $\eta$. Moreover, the classical Bers' theorem can be derived from the assumption that every hyperbolic surface admits a proper essential subsurface with definite entropy.  
\end{theorem}

Given a graph $G$, let $\rk \pi_1(G)$ denote the rank of its fundamental group. Our main theorem is the direct analog of the first sentence of Theorem \ref{thm:BersForGraphs} in the setting of metric graphs:

\begin{theorem}\label{thm:main}
For each $r \geq 3$, there is a constant $C_r > 0$ so that given a graph $G$ with $\rk \pi_1(G) = r \geq 3$ and a length function $\ell$ with unit entropy there is a proper subgraph $G' \subset G$ such that the restriction of $\ell$ to $G'$ has entropy at least $C_r$.
\end{theorem}

In addition to achieving a Bers Lemma like result for metric graphs, we place Theorem \ref{thm:main} in conversation with the work Kim-Lim \cite{Kim2020}, which establishes a relation between the entropy of a graph $G'$ obtained from a graph $G$ by adding an additional edge between a pair of initially non-adjacent vertices. In some sense, our work concerns the opposite direction of starting with a unit entropy graph, deleting an edge, and estimating (from below) the entropy of the resulting graph.

As mentioned above, in~\cite{ar:ACR23}, the authors study a piecewise Riemannian metric on the outer space of a free group whose definition is inspired by the classical Weil--Petersson metric on the Teichm{\"u}ller space. We show that, as in the case of the Weil--Petersson metric, this metric is incomplete. Unlike the Weil--Petersson metric and the associated mapping class group action, we determine that so long as the rank is at least $4$, the completion admits a global fixed point for the $\Out(F_{r})$ action. In fact, we conjectured that the entire outer space has finite diameter, but in attempting to prove this conjecture we ran into the need of a lemma along the lines of Theorem \ref{thm:main}.  Indeed, points in the completion correspond to unit entropy metrics on proper subgraphs and measuring the distance to such a point is related to how much scaling is needed to increase the entropy on the subgraph to 1.  Knowing that there is a subgraph whose entropy is bounded away from 0 gives a candidate for a point in the completion at a bounded distance.  However, there is still some analysis necessary to bound the integral that computes this distance.  We therefore reiterate the following conjecture: 

\begin{conjecture} The entropy metric (see \cite{ar:ACR23}) has finite diameter on the outer space for all ranks $r \ge 4$.
\end{conjecture}

We note that, in Section \ref{sec:roses} where we prove Theorem \ref{thm:main} for roses, we in fact show that the constant $C_{r}$ for roses can actually be taken to converge to $1$ as $r \rightarrow \infty$. We conjecture that this is true for general metric graphs: 

\begin{conjecture} The constant $C_{r}$ from Theorem \ref{thm:main} is universally bounded from below in $r$, and moreover converges to $1$ as $r \rightarrow \infty$. 
\end{conjecture} 

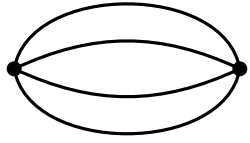
\begin{figure}[ht]
\begin{tikzpicture}[scale=1.5]
\filldraw (-1.5,0) circle [radius=1.75pt];
\filldraw (0.5,0) circle [radius=1.75pt];
\draw[very thick] (-1.5,0) to[out=25,in=155] (0.5,0);
\draw[very thick] (-1.5,0) to[out=-25,in=-155] (0.5,0);
\draw[very thick] (-1.5,0) to[out=80,in=100] (0.5,0);
\draw[very thick] (-1.5,0) to[out=-80,in=-100] (0.5,0);
\end{tikzpicture}
\caption{The graph $G$}\label{fig:conjecture graph}
\end{figure}

In fact, we conjecture that the most extreme case occurs when $G$ is the rank 3 graph with two vertices and four edges as shown in Figure~\ref{fig:conjecture graph} and $\ell = (\log(3),\log(3),\log(3),\log(3))$.  Removing any of the four edges results in a proper subgraph $G'$ with entropy $\log(2)/\log(3) \approx 0.6309...$.

\subsection*{Outline of Paper} In Section \ref{sec:prelim}, we recall the basic properties of entropy that we will need, and we do some computations to obtain some helpful inequalities relating the entropies of various roses and barbells. We also spell out a proof of Theorem \ref{thm:BersForGraphs}. In Section \ref{sec:estimate}, we define a function encoding the scalar one should multiply against the lengths of all but one edge $e$ of a graph in order to maintain entropy $1$ while blowing up the length of $e$, and we relate the entropy of the subgraph $G - e$ to the derivative of this function. In Section \ref{sec:roses}, we prove Theorem \ref{thm:main} for roses.  In Section \ref{sec:Reduction}, we outline a reduction in the argument for the general proof of Theorem \ref{thm:main} which will allow us to assume that every non-loop edge of $G$ is incident to a self-loop on both sides. In Section \ref{sec:lemma}, we prove an upper bound on the equilibrium measure of a non-loop edge $e$ in terms of the measures associated to the self-loops incident to either side of $e$. Finally, in Section \ref{sec:proof}, we use the bounds in Section \ref{sec:lemma} complete the proof of Theorem \ref{thm:main}.

\subsection*{Acknowledgements} The first author would like to thank Samuel Taylor for a number of very helpful conversations about this problem over the course of the last two years. He would also like to thank former students Atira Glenn-Keough, Peter Kulowiec, Halley Kucirka, and Rachel Niebler who focused on this and related problems for their senior thesis work during the 2023-2024 academic year. The second author would like to thank Yo'av Rieck, Chris Cashen, and Derrick Wigglesworth for many helpful conversations.

\section{Preliminaries}\label{sec:prelim}

In this section we collect the necessary definitions and facts needed for the sequel.  Much of this material is taken from the works of Aougab--Clay--Rieck~\cite{ar:ACR23}, Parry--Pollicott~\cite{ar:PP90}, and Pollicott--Sharp~\cite{ar:PS14}.  We refer to those sources for motivation and further discussion on these topics.  
In addition, we give a computation in Section~\ref{subsec:computations} that will be used in the proof of Theorem~\ref{thm:main}.

\subsection{Graphs}\label{subsec:graphs}

A \emph{graph} is a tuple $G = (V,E,o,\tau,\bar{\phantom{e}})$ where:
\begin{enumerate}
    \item $V$ and $E$ are sets, called the \emph{vertices} and the \emph{directed edges}, 
\item  $o,\tau \from E \to V$ are functions that specify the \emph{originating} and \emph{terminating} vertices of an edge, and
\item $\bar{\phantom{e}} \from E \to E$ is a fixed-point free involution such that $o(e) = \tau(\bar{e})$.
\end{enumerate}
Given a graph, we fix a subset $E_+ \subset E$ that consists of exactly one edge from each pair $\{e,\overline{e}\}$.

A \emph{length function} on $G$ is a function $\ell \from E_+ \to \RR_{>0}$.  The moduli space of all such functions is denoted $\calM(G)$.  A length function $\ell$ extends to a function on $E$ by defining $\ell(e) = \ell(\bar{e})$ if $e \notin E_+$.

A \emph{\textup{(}based\textup{)} circuit} is a finite sequence of edges $(e_1,\ldots,e_n)$ where:
\begin{enumerate}
\item $t(e_i) = o(e_{i+1})$ for $i = 1, \ldots, n-1$ and $t(e_n) = o(e_1)$, and
\item $e_i \neq \bar{e}_{i+1}$ for $i = 1, \ldots, n-1$ and $e_n \neq \bar{e}_1$.
\end{enumerate}
Given a length function $\ell \in \calM(G)$, the length of a circuit $\gamma = (e_1,\ldots,e_n)$ is defined as:
\begin{equation*}
\ell(\gamma) = \sum_{i=1}^n \ell(e_i).
\end{equation*}
The set of all circuits in $G$ is denoted by $\calC(G)$.  For a given length function $\ell \in \calM(G)$, the subset of circuits whose length is at most $t$ is denoted $\calC_{G,\ell}(t)$.

\subsection{Entropy}\label{subsec:entropy}

The \emph{entropy} of a length function $\ell \in \calM(G)$ is defined as:
\begin{equation*}
\entropy_{G}(\ell) = \lim_{t \rightarrow \infty} \frac{1}{t} \log \abs{\calC_{G,\ell}(t)}.
\end{equation*}
This function $\entropy_{G} \from \calM(G) \to \RR_{\geq 0}$ is real analytic and strictly convex~\cite[Proposition~A.4]{ar:McMullen15}.  It is readily verified that $\entropy_G$ is homogeneous of degree $-1$.  That is:

\begin{equation} \label{homogeneous}
\entropy_G(\alpha\ell) = \frac{1}{\alpha}\entropy_G(\ell), 
\end{equation}
for any $\alpha > 0$~\cite[Lemma~3.4]{ar:ACR23}.

The subset of length functions with unit entropy is denoted $\calM^1(G)$.  That is:
\begin{equation*}
\calM^1(G) = \{ \ell \in \calM(G) \mid \entropy_G(\ell) = 1\}.
\end{equation*}

We associate to the graph $G$ a matrix that records the incident relations.  Specifically, we define an $\abs{E} \times \abs{E}$ matrix $A_{G}$ by
\begin{equation*}
A_{G}(e,e') = \begin{cases} 1 & \mbox{ if } \tau(e) = o(e') \mbox{ and } \bar{e} \neq e', \\
0 & \mbox{ else. }
\end{cases}
\end{equation*}
Given a length function $\ell \in \calM(G)$, we set $A_{G,\ell}$ to be the matrix obtained by multiplying the row of $A_G$ corresponding to $e \in E$ by $\exp(-\ell(e))$ for each edge.  That is:
\begin{equation*}
A_{G,\ell}(e,e') = \exp(-\ell(e))A_G(e,e').
\end{equation*}
This matrix is irreducible.  Moreover, we have that $\spec(A_{G,h\ell}) = 1$ precisely when $h = \entropy_G(\ell)$ where $\spec(\param)$ denotes the spectral radius (cf.~\cite[Lemma~3.1(2)]{ar:PS14}).  In particular:  
\begin{equation}\label{eq:entropy 1}
\entropy_G(\ell) = 1 \Leftrightarrow \spec(A_{G,\ell}) = 1
\end{equation}

It follows from \ref{eq:entropy 1} and the Perron-Frobenius theorem that when $\entropy_G(\ell) = 1$, there are positive vectors $\bu,\bv \in \RR^{2\abs{E}}$ such that:
\begin{equation}\label{eq:PF vectors}
\bu^TA_{G,\ell} = \bu^T, \quad A_{G,\ell}\bv = \bv, \text{ and } \bu^T\bv = 1. 
\end{equation}
Moreover, this pair of vectors is unique up to scale by $\alpha, \frac{1}{\alpha}$.  The \emph{equilibrium measure} $\mu$ is defined by taking the entry-wise products of $\bu$ and $\bv$, i.e., $\mu(e) = \bu(e)\bv(e)$.  This determines a measure on the set of bi-infinite lines in $G$ in the usual way, c.f.~\cite[Section~6.6]{bk:Walters82}.  We remark that, by symmetry, one has $\mu(e) = \mu(\bar{e})$.

We define a function $F_G \from \calM(G) \to \RR$ by:
\begin{equation}\label{eq:hypersurface}
F_G(\ell) = \det\left(I - A_{G,\ell}\right).
\end{equation}
From~\eqref{eq:entropy 1}, it follows that $F_G(\ell) = 0$ if $\entropy_G(\ell) = 1$.  Furthermore, we have the following relation between $F_G$ and the equilibrium measure $\mu$. 

\begin{lemma}\label{lem:proportion}
Let $\ell \in \calM^1(G)$ and let $\mu$ be the corresponding equilibrium measure.  Then $\nabla F(\ell)$ is proportional to $\mu$.  That is, there is a nonzero constant $C$ such that $\mu(e) = C\partial_eF(\ell)$ for all $e \in E_+$.
\end{lemma}

\begin{proof}
This follows as $\calM^1(G)$ is a component of $F_G^{-1}(0)$ and also the zero locus of $\log \spec(A_{G,\ell})$ since the gradient of the latter function is $\mu$.  See the proofs of Lemmas 3.9 and 4.4 in~\cite{ar:ACR23} for complete details.
\end{proof}

Finally, we note that $\entropy_G$ extends to a continuous function where we allow $\infty$ as a value for the length of an edge.  Indeed, this follows from the implicit function theorem as entropy is characterized by $\spec(A_{G,\entropy_G(\ell)\ell}) = 1$, and since characteristic functions and the Perron--Froebinous eigenvalue depend continuously on the matrix entries.

\subsection{Connection to Bers' constant in hyperbolic geoemetry}

Before using the terms and tools introduced above to carry out some computations that will be needed to prove the main theorem, we give a poof of Theorem \ref{thm:BersForGraphs} from the introduction: 

\begin{theorem} Given $\chi < 0$ and $L \ge 0$, there is a constant $\eta >0$ such that if $X$ is a hyperbolic surface (potentially with totally geodesic boundary) with Euler characteristic $\chi$ and which is not a pair of pants and with all boundary components of length at most $L$, there is a proper essential subsurface $Y \subsetneq X$ with entropy at least $\eta$. Conversely, if one assumes that every such surface admits a proper subsurface with entropy at least $\eta = \eta(\chi, L)$, then one can prove a version of the classical Bers' theorem, i.e., that there exists a pants decomposition with total length bounded above in terms of only $\chi,L$, and $\eta$. 
\end{theorem}

\begin{proof} By Bers' theorem, there is a pants decomposition on $X$ with total length at most $\mathcal{L} = \mathcal{L}(\chi, L)$, for some constant depending only on the Euler characteristic $\chi$ and on the length of the boundary, $L$. Focus on one pair of pants $P$ in the complement of this pants decomposition: it is a hyperbolic surface in its own right, with three totally geodesic boundary components each with length at most $\mathcal{L}$. The entropy $h(P)$ can be expressed as 
\[ h(P) = \lim_{R \rightarrow \infty} \frac{1}{R} \cdot \#\left\{ \gamma \in \Gamma = \langle A, B \rangle: d_{\mathbb{H}}(o, \gamma \cdot o) \leq R \right\}, \]
where $\Gamma < PSL_{2}(\mathbb{R})$ is the given image of the representation of $\pi_{1}(P) \cong F_{2}$; $A,B$ are hyperbolic isometries corresponding to two of three boundary components of $P$; and $o \in \mathbb{H}$ is a chosen base-point. 

Letting $|\gamma|$ denote the word length of $\gamma$ in the basis $A,B$, by a straightforward application of the triangle inequality one has that 
\[ d_{\mathbb{H}}(o, \gamma \cdot o) \leq \mathcal{L} \cdot |\gamma|. \]
Thus, 
\[ \left\{\gamma \in \Gamma: |\gamma| \leq R/\mathcal{L} \right\} \subset \left\{ \gamma \in \Gamma: d_{\mathbb{H}}(o, \gamma \cdot o) \leq R \right\}.\]
Therefore, the number of elements in $\pi_{1}(P)$ with word length at most $R/\mathcal{L}$ is a lower bound for the number of elements that displace $o$ by at most $R$. By a standard count, the number of words in $F_{2}$ with word length at most $n$ is on the order of $3^{n}$, so we conclude that 
\[ 3^{R/\mathcal{L}} \leq \# \left\{ \gamma \in \Gamma: d_{\mathbb{H}}(o, \gamma \cdot o) \leq R \right\}. \]
Taking logs of both sides and then letting $R \rightarrow \infty$ yields a lower bound on $h(P)$ depending only on $\mathcal{L}$, as desired. 

Conversely, assume that every hyperbolic surface has a proper subsurface with entropy at least some $\eta >0$ depending only on its Euler characteristic $\chi$ and on the length of its boundary, $L$. The standard proof of Bers' theorem (see for instance 12.4.2 of \cite{Farb2012}) uses the Gauss-Bonnet theorem and its consequence that $\mbox{Area}(X)$ is a function of $\chi(X)$ and $\ell(\partial X)$. Our argument will follow a similar structure, but within that structure we will replace the Gauss Bonnet theorem and considerations of area with the assumption that $X$ admits a proper subsurface with definite entropy $\eta$. 

Now, let $Y \subsetneq X$ be the proper subsurface with entropy $\eta$. It follows that $X$ itself has entropy at least $\eta$. For some constant $M >0$ assume that every simple closed geodesic (and therefore every closed geodesic) in $X$ has length at least $M$. For $R$ very large, consider the ball $B(o, R)$ of radius $R$ in the universal cover, $\tilde{X}$, about $o$ (note that $\tilde{X}$ needn't be equal to all of $\mathbb{H}^{2}$). One can choose $R$ sufficiently large so that the logarithm of the number of orbit points corresponding to the action of $\pi_{1}(X)$ in $B(o, R)$, divided by $R$, is within some small $\epsilon$ of $h(X) > \eta$. We remark that this value of $R$ depends only on $\chi(X)$ and not on the particular choice of metric; indeed, because $h$ varies analytically over the moduli space and the $M$-thick part is compact, the rate of convergence of the limit in the definition of $h$ is bounded from below in terms only of $\chi$.

Now, for any two such orbit points $x,y$, one has that $B(x,M) \cap B(y,M) = \emptyset$. In the event that $x$ is within $M$ of exactly one boundary component of $X$, then $B(x,M)$ contributes an area of at least $\mbox{Area}(B(x,M))/2$. Suppose first that no $x$ is within $M$ of more than one such boundary component. Then the total number of orbit points is at most 
 
\[ \frac{\mbox{Area}(B(o,R))}{2 \cdot \mbox{Area}(B(o,M))}. \]
It follows that 
\[ \exp(h \cdot R) \le \frac{\cosh(R)-1}{2(\cosh(M)-1)} \sim \exp(R-M), \]
which in turn implies an upper bound for $M$ in terms of $R$ and $h$ (and therefore in terms of $\chi$ and $h$), as desired. 

On the other hand, when some $x$ is within $M$ of two distinct lifts of boundary components of $X$, then there is a simple closed geodesic on $X$ with length at most $4M + 2 \ell(\partial X)$: one simply travels from $x$ to one boundary component, traverses it, then returns to $x$ and then moves to the second boundary component, traverses it, and then returns to $x$ to close up. 

Thus, we have succeeded in producing a closed geodesic $\alpha$ on $X$ with length bounded above in terms of only $\chi, L$, and $\eta$. The existence of the desired pants decomposition follows by a standard inductive argument in which one cuts along $\alpha$ and then iterates. 
\end{proof}

\subsection{Computations with entropy}\label{subsec:computations}

We end this section with calculations and bounds on the entropy function that will be needed in the proof of Theorem~\ref{thm:main}.

\subsubsection{The 2--Rose} \label{subsec:rose}

\begin{figure}[ht]
\centering
\begin{tikzpicture}[scale=1.5]
\filldraw (-1.5,0) circle [radius=1.75pt];
\draw[very thick] (-2.15,0) circle[radius=0.65cm];
\draw[very thick] (-0.85,0) circle[radius=0.65cm];
\node at (-2.15,-0.85) {\footnotesize $e_1$};
\node at (-0.85,-0.85) {\footnotesize $e_2$};
\node at (-1.7,0) {\footnotesize $v$};
\end{tikzpicture}
\caption{The 2--rose graph $\calR_2$.}\label{fig:2-rose}
\end{figure}
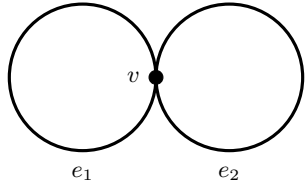

We denote the edges of the 2--rose $\calR_2$ by $e_1$ and $e_2$, and the vertex by $v$.  See Figure~\ref{fig:2-rose}.  

The moduli space $\calM^1(\calR_2) \subset \RR^2_{>0}$ is the zero set of the function:
\begin{equation*}
F_{\calR_2}(a,b) =1-  \exp(-a) - \exp(-b) - 3\exp(-a-b), 
\end{equation*}
where $a,b$ denote the lengths of $e_{1}, e_{2}$ respectively. 
(cf.~\cite[Section~6.1]{ar:PS14} or \cite[Section~6.1]{ar:ACR23})  In particular, if we define
\begin{equation}\label{eq:R2}
R_2(x) = -\log\left(\frac{1-\exp(-x)}{1+3\exp(-x)}\right)
\end{equation}
then $\entropy_{\calR_2}(x,R_2(x)) = 1$.

\subsubsection{The Barbell} \label{subsec:barbell}

\begin{figure}[ht]
\centering
\begin{tikzpicture}[scale=1.5]
\draw[very thick] (-1.5,0) -- (-0.5,0);
\filldraw (-1.5,0) circle [radius=1.75pt];
\filldraw (-0.5,0) circle [radius=1.75pt];
\draw[very thick] (-2.15,0) circle[radius=0.65cm];
\draw[very thick] (0.15,0) circle[radius=0.65cm];
\node at (-2.15,-0.85) {\footnotesize $e_1$};
\node at (0.25,-0.85) {\footnotesize $e_2$};
\node at (-1,-0.2) {\footnotesize $e_3$};
\node at (-1.7,0) {\footnotesize $v$};
\node at (-0.3,0) {\footnotesize $w$};
\end{tikzpicture}
\caption{The barbell graph $\calB_2$.}\label{fig:barbell}
\end{figure}
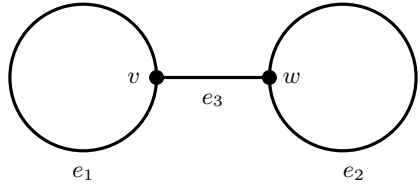

The graph $\calB_2$ has two vertices $v$ and $w$ and three edges $e_1$, $e_2$, and $e_3$ where $o(e_1) = t(e_1) = v$,  $o(e_2) = t(e_2) = v$, $o(e_3) = v$, and $t(e_3) = w$.  See Figure~\ref{fig:barbell}.  

We will provide a condition on the length of the edges in $\calB_2$ that ensures entropy of at least $\frac{1}{5}$.  The following lemma will allow us to compare the entropy of a given barbell to that of a certain rose: 

\begin{lemma}\label{lem:rose barbell compare}
For any positive constants, $a$, $b$, and $c$, we have that $\entropy_{\calB_2}(a,b,c) \geq \entropy_{\calR_2}(a,b+2c)$.
\end{lemma}

\begin{proof}
Fix a length function $\ell = (a,b,c) \in \calM(\calB_2)$ and let $\ell' \in \calM(\calR_2)$ be the length function defined by $\ell' = (a, b+2c)$.  

By identifying the initial and terminal segments of length $c$ of the edge $e_2$ in $\calR_2$, we get a map $\calR_2 \to \calB_2$ that induces a length decreasing bijection on cycles $\calC(\calR_2) \to \calC(\calB_2)$.  Hence we find that $\abs{\calC_{\calR_2,\ell'}(t)} \leq \abs{\calC_{\calB_2,\ell}(t)}$ and from this the lemma follows from the definition of entropy.
\end{proof}

Using this comparison, we can now prove a lemma that will be used in the proof of Theorem~\ref{thm:main} to produce a subgraph of definite entropy. 

\begin{lemma}\label{lem:barbell estimate}
Let $\ell \in \calM(\calB_2)$ be a length function where $\ell(e_1) = 3\exp(-c/2)$, $\ell(e_2) = c/4$ and $\ell(e_3) = c$ for any $c > 0$.  Then $\entropy_{\calB_2}(\ell) \geq \frac{1}{5}$.
\end{lemma}

\begin{proof}
Fix a positive constant $c$ and let $\ell \in \calM(\calB_2)$ be as in the statement.  Let $a = 3\exp(-c/2)$.  Then $\ell(e_2) = -\frac{1}{2}\log(a/3)$ and $\ell(e_3) = -2\log(a/3)$.  By Lemma~\ref{lem:rose barbell compare}, we find that:
\begin{equation}\label{eq:rose barbell comparison}
\entropy_{\calB_2}(\ell) = \entropy_{\calB_2}(a,-\frac{1}{2}\log(a/3), -2\log(a/3)) \geq \entropy_{\calR_2}(a,-5\log(a/3)).
\end{equation}  

Let $R_2(a) = -\log\left(\frac{1-\exp(-a)}{1+3\exp(-a)}\right)$ as in~\eqref{eq:R2}.  Then $\entropy_{\calR_2}(a,R_2(a)) = 1$ as explained in Section~\ref{subsec:rose}. 

\begin{claim*}
For $0 < a < 3$, we have $-5\log(a/3) \leq 5R_2(a)$.
\end{claim*}

\begin{proof}[Proof of Claim]
This is equivalent to
\begin{equation}\label{eq:inequality}
\exp(a)(1-a/3) \leq 1+a.
\end{equation}  
The inequality~\eqref{eq:inequality} is true for $a = 0$.  The derivative of the left hand side of~\eqref{eq:inequality} is $\frac{1}{3}(2-a)\exp(a)$.  We observe that 
\begin{equation}\label{eq:inequality 2}
\frac{1}{3}(2-a)\exp(a) \leq 1
\end{equation} 
for $0 \leq a \leq 3$, hence proving the claim, as 1 is the derivative of the right hand side of~\eqref{eq:inequality}.  Indeed, the function $\frac{1}{3}(2-a)\exp(a)$ on $0 \leq a \leq 3$ is maximized at $a= 1$ and we observe that $\frac{1}{3}\exp(1) \leq 1$ holds. 
\end{proof}

Therefore, by the claim, we find that:
\begin{equation*}
\entropy_{\calR_2}(a,-5\log(a/3)) \geq \entropy_{\calR_2}(a,5R_2(a)) \geq \entropy_{\calR_2}(5a,5R_2(a)) = \frac{1}{5}.
\end{equation*}  
The lemma now follows by combining this with~\eqref{eq:rose barbell comparison}.
\end{proof}

\section{Computing the Entropy of a Subgraph} \label{sec:estimate}

Let $e$ be any edge of a graph $G$ and fix a length function $\ell \in \calM^1(G)$.  We will be interested in a smooth family of length functions  $\psi_{t} \in \calM^1(G)$ defined as follows for $0 \leq t < \infty$:  
\begin{equation}\label{eq:psi_t}
\psi_t(e') = \begin{cases}
\ell(e) + t & \text{ if } e' = e \\
j(t) \cdot \ell(e') & \text{ otherwise}.
\end{cases}
\end{equation} 
where $j\from \RR_{\geq 0} \to \RR_{\geq 0}$ is defined to be the function so that $\entropy_{G}(\psi_t) = 1$. In other words, we linearly increase the length of $e$ and scale the remaining edges to maintain unit entropy.  The smoothness of $j$ follows from the smoothness of the entropy function (for instance, by the implicit function theorem).  We call the path $\psi_t \in \calM^1(G)$ a \emph{linear time blow up of $\ell$ along $e$}. 

The main result of this section shows that we can compute the entropy $\ell$ restricted to the subgraph $G' = G - e$ using the derivative $j'(t)$.  Further, we compute a formula for $j'(t)$ that we use in the later sections.

\begin{proposition}\label{prop:estimate}
Let $\ell \in \calM^1(G)$ where $G$ is a connected graph with $\rk \pi_1 G \geq 3$, and fix an edge $e \in E(G)$.  Let $\psi_t \in \calM^1(G)$ be the linear time blow-up of $\ell$ along $e$ with corresponding scaling function $j \from \RR_{\geq 0} \to \RR_{\geq 0}$, and let $\mu_t$ denote the equilibrium measure for the length function $\psi_t$.
\begin{enumerate}
\item \label{item:derivative}The derivative of $j$ satisfies:
\begin{equation}  \label{eq:derivative of j}
j'(t) = \frac{- \mu_t(e) }{\sum_{e' \neq e} \ell(e')\mu_t(e')}.
\end{equation}
\item \label{item:entropy} For the subgraph $G' = G - e$ we have:
\begin{equation}\label{eq:entropy estimate}
\entropy_{G'}(\ell\big|_{G'}) = 1 - \int_{0}^{\infty} \abs{j'(t)} dt. 
\end{equation}
\end{enumerate}
\end{proposition}

\begin{proof}
First we show~\eqref{item:derivative}. We can compute an expression for $j'(t)$ from the equation $F_{G}(\psi_{t}) = 0$.  Indeed, differentiating both sides with respect to $t$ yields 
\begin{equation*} \sum_{e' \in E(G)} \frac{d}{dt} \psi_t(e') \cdot \partial_{e'} F(\psi_{t}) = 0.
\end{equation*}
By~\eqref{eq:psi_t} we have that:
\begin{equation*}
\frac{d}{dt}\psi_t(e') = \begin{cases}
1 & \text{ if }e' = e, \\
j'(t) \cdot \ell(e') & \text{ otherwise}.
\end{cases}
\end{equation*}
Therefore, we obtain the equation 
\begin{equation*} 
- \partial_{e} F_G(\psi_{t}) = \sum_{e' \neq e} j'(t) \cdot \ell(e') \partial_{e'} F_G(\psi_{t}).
\end{equation*}
Solving for $j'(t)$ yields 
\begin{equation*}
j'(t) = \frac{ - \partial_{e} F_G(\psi_{t})}{\sum_{e' \neq e} \ell(e') \partial_{e'} F_G(\psi_{t})}.
\end{equation*}
Finally, the fact that $\nabla F_G(\psi_t)$ is proportional to $\mu_{t}$, the equilibrium measure at $\psi_t$, (Lemma~\ref{lem:proportion}) yields the formula in \eqref{eq:derivative of j}.

Next we demonstrate~\eqref{item:entropy}.   We observe that $j(t)$ is monotonically decreasing and bounded below by $0$.  This is clear from the definitions as we must shrink the edges in $G' = G - e$ more as the length of $e$ grows, but the explicit formula in~\eqref{eq:derivative of j} shows that $j'(t) < 0$ as well.  This implies that $\lim_{t \to \infty} j(t)$ exists.  Denote this limit by $j_\infty$.  We also observe that $j(0) = 1$ as $\entropy_G(\ell) = 1$.    

Hence, for the right hand side of \eqref{eq:entropy estimate}, we find:
\begin{equation*} 
1 - \int_{0}^{\infty} \abs{j'(t)} dt = 1 + \int_{0}^{\infty} j'(t) dt = 1 + \lim_{t \to \infty} j(t) - j(0) = j_\infty.
\end{equation*}

As the rank of $G$ is at least 3, $j_\infty > 0$.  Indeed, at least one component of $G'$ has rank at least 2 and hence $\entropy_{G'}\left(\ell\big|_{G'}\right) > 0$.  This forces $j_t > 1/\entropy_{G'}\left(\ell\big|_{G'}\right)$ for all $t$.

By continuity of the extended entropy function and the definition of $\psi_t$, we find:
\begin{equation*}
1 = \lim_{t \to \infty} \entropy_G(\psi_t) = \entropy_{G'}(j_\infty \ell\big|_{G'}) = \frac{1}{j_\infty}\entropy_{G'}(\ell\big|_{G'}).
\end{equation*}
Thus $\entropy_{G'}(\ell\big|_{G'}) = j_\infty$, which verifies \eqref{eq:entropy estimate}.
\end{proof}

Therefore, our approach for proving Theorem \ref{thm:main} will be to bound the improper integral $\displaystyle \int_0^\infty \abs{j'(t)} \, dt$ away from 1.

\section{Roses}\label{sec:roses}

In this section we will prove Theorem~\ref{thm:main} for $r$--roses $\calR_r$, $r \geq 3$.  As stated in the introduction, our method of proof for the general case requires a non-loop edge.  Hence we must deal with the rose case separately.  The strategy here though is similar to the general case and provides a good warm-up and overview for the methods used later on.  Namely, we will make use of the structure of the matrix $A_{\calR_r}$ to get bounds on the values of the equilibrium measure related to edge lengths that, via Proposition~\ref{prop:estimate}, allow us to estimate the entropy of a subgraph.

These bounds are the context of the following lemma.

\begin{lemma}\label{lem:rose estimate}
Fix $\ell \in \calM^1(\calR_r)$ with equilibrium measure $\mu$ and fix two edges $e_1,e_2 \in E(\calR_r)$.  Then:
\begin{equation*}
\exp(\ell(e_1))\mu(e_1) < 4\exp(\ell(e_2))\mu(e_2)
\end{equation*}
\end{lemma}

\begin{proof}
Let $A = A_{\calR_r,\ell}$ and fix vectors $\bu$ and $\bv$ such that $\bu^T A = \bu^T$, $A\bv = \bv$, and $\bu^T\bv = 1$.  Then we have that $\mu(e) = \bu(e) \bv(e)$ for all edges as explained in Section~\ref{subsec:entropy}.   Throughout the calculations below we will make use of the fact that $\bu(e_i) = \bu(\bar{e}_i)$ for any edge $e_i$ and similarly for $\bv$ and $\mu$.

As $\bu^TA = \bu^T$, we find that:
\begin{align*}
\bu(e_1) & = \exp(-\ell(e_1))\bu(e_1) + 2\sum_{i\neq 1} \exp(-\ell(e_i))\bu(e_i), \text{and} \\
\bu(e_2) & = \exp(-\ell(e_2))\bu(e_2) + 2\sum_{i\neq 2} \exp(-\ell(e_i))\bu(e_i).
\end{align*}
From this we observe that:
\begin{equation}\label{eq:u1 < u2 rose}
\bu(e_1) < 2\bu(e_2).
\end{equation}

As $A\bv= \bv$, we find that:
\begin{align*}
\bv(e_1) & = \exp(-\ell(e_1))\left(\bv(e_1) + 2\sum_{i\neq 1} \bv(e_i)\right), \text{ and} \\
\bv(e_2) & = \exp(-\ell(e_2))\left(\bv(e_2) + 2\sum_{i\neq 2} \bv(e_i)\right). \\
\end{align*} 
As above, this implies that:
\begin{equation}\label{eq:v1 < v2 rose}
\exp(\ell(e_1))\bv(e_1) < 2\exp(\ell(e_2))\bv(e_2).
\end{equation}

Combining \eqref{eq:u1 < u2 rose} and \eqref{eq:v1 < v2 rose} we find:
\begin{equation*}
\exp(\ell(e_1))\mu(e_1) = \exp(\ell(e_1))\bu(e_1)\bv(e_1) < 4\exp(\ell(e_2))\bu(e_2)\bv(e_2) = 4\exp(\ell(e_2))\mu(e_2). \qedhere
\end{equation*}
\end{proof}

Next, we will show how these bounds prove Theorem~\ref{thm:main} in the case of roses.

\begin{proposition}\label{prop:rose}
For each $r \geq 3$, there is a constant $B_r > 0$ so that for any given length function $\ell \in \calM^1(\calR_r)$, there is a proper subgraph $G' \subset \calR_r$ such that $\entropy_{G'}\left(\ell\big|_{G'}\right)\geq B_r$.  
Moreover, $\lim_{r \to \infty} B_r = 1$.
\end{proposition}

\begin{proof}
Fix a length function $\ell \in \calM^1(\calR_r)$ and rename the edges so that $\ell(e_1) \geq \ell(e_2) \geq \ell(e_i)$ for $3 \leq i \leq r$.  Let $\psi_t \in \calM^1(\calR_r)$ be the linear time blow-up of $\ell$ along $e_1$ with corresponding scaling function $j(t)$.  Let $\mu_t$ be the equilibrium measure at $\psi_t$ so that we are in set-up of Proposition~\ref{prop:estimate}.

By Proposition~\ref{prop:estimate}\eqref{item:derivative}, we have:
\begin{equation*}
\int_0^\infty \abs{j'(t)} \, dt = \int_0^\infty \frac{\mu_t(e_1)}{\sum_{i\neq 1}\ell(e_i)\mu_t(e_i)} \, dt 
\end{equation*}
Lemma~\ref{lem:rose estimate} gives us that $\mu_t(e_1) < 4\exp(j(t)\ell(e_2) - \ell(e_1)-t)\mu_t(e_2)$, and therefore:
\begin{equation*}
\int_0^\infty \frac{\mu_t(e_1)}{\sum_{i\neq 1}\mu_t(e_i)\ell(e_i)} \, dt < 4\int_0^\infty \frac{\exp(j(t)\ell(e_2)-\ell(e_1)-t)\mu_t(e_2)}{\sum_{i\neq 1}\ell(e_i)\mu_t(e_i)} \, dt 
\end{equation*}
As $j(t)\ell(e_2) \leq \ell(e_2) \leq \ell(e_1)$, we have $\exp(j(t)\ell(e_2)-\ell(e_1)) \leq 1$.  This gives:
\begin{equation*}
4\int_0^\infty \frac{\exp(j(t)\ell(e_2)-\ell(e_1)-t)\mu_t(e_2)}{\sum_{i\neq 1}\ell(e_i)\mu_t(e_i)} \, dt \leq 4\int_0^\infty \exp(-t) \frac{\mu_t(e_2)}{\sum_{i\neq 1}\ell(e_i)\mu_t(e_i)} \, dt
\end{equation*}
Lastly, since $\ell(e_2)\mu_t(e_2) \leq \sum_{i\neq 1}\ell(e_i)\mu_t(e_i)$, we conclude:
\begin{align*}
\int_0^\infty \abs{j'(t)} \, dt & \leq  4\int_0^\infty \exp(-t) \frac{\mu_t(e_2)}{\sum_{i\neq 1}\ell(e_i)\mu_t(e_i)} \, dt \\
& =  4\int_0^\infty \frac{\exp(-t)}{\ell(e_2)} \frac{\mu_t(e_2)\ell(e_2)}{\sum_{i\neq 1}\ell(e_i)\mu_t(e_i)} \, dt \\
& \leq \frac{4}{\ell(e_2)} \int_0^\infty \exp(-t) \, dt \\
& = \frac{4}{\ell(e_2)}. 
\end{align*}

Let $G' = \calR_r - e_1$.  Hence by Proposition~\ref{prop:estimate}\eqref{item:entropy} we find:
\begin{equation}\label{eq:edge bound}
\entropy_{G'}(\ell\big|_{G'}) = 1 - \int_0^\infty \abs{j'(t)} \, dt \geq 1 - \frac{4}{\ell(e_2)}.
\end{equation}

Notice that we must have $\ell(e_2) \geq \log(2r-3)$.  Indeed, the length function $\ell$ restricted to the subrose $\calR_r - e_1$ has entropy less than 1.  For $r \geq 29$, we have $\log(2r-3) > 4$ and therefore, \eqref{eq:edge bound} gives:
\begin{equation*}
\entropy_{G'}(\ell\big|_{G'}) \geq 1 - \frac{4}{\ell(e_2)} \geq 1 - \frac{4}{\log(2r-3)} > 0.
\end{equation*}

When $r < 29$, we will need to consider two subcases.  First, if $\ell(e_2) > 5$, then \eqref{eq:edge bound} gives:
\begin{equation*}
\entropy_{G'}(\ell\big|_{G'}) \geq 1 - \frac{4}{5} = \frac{1}{5}.
\end{equation*}
Next, we suppose that $\ell(e_2) \leq 5$.  Thus $\ell(e_i) \leq 5$ for $2 \leq i \leq r$.  In this case, we find that:
\begin{align*}
\entropy_{G'}(\ell\big|_{G'}) & \geq \entropy_{\calR_{r-1}}(5,\ldots,5) \\ 
&= \frac{\log(2r-3)}{5}\entropy_{\calR_{r-1}}(\log(2r-3),\ldots,\log(2r-3)) \\
& = \frac{\log(2r-3)}{5} \\
& > \frac{1}{5}.
\end{align*}

Hence we can take $B_r = \frac{1}{5}$ for $3 \leq r < 29$, and $B_r = 1 - \frac{4}{\log(2r-3)}$ for $r \geq 29$.
\end{proof}

We conjecture that $B_3 = \log(3)/\log(5) \approx 0.6826...$ which is obtained for the uniform length function $\ell = (\log(5),\log(5),\log(5)) \in \calM^1(\calR_3)$.  Further, we conjecture that $B_r \geq B_3$ for all $r \geq 3$.

\section{Reduction} \label{sec:Reduction}

\begin{notation}\label{notation:entropy}
We define the following quantities:
\begin{align*}
\entropysup_G(\ell) & = \max\left\{ \entropy_{G'}\left(\ell\big|_{G'}\right) \mid G' \subset G \text{ is a proper subgraph}\right\} \\
\entropyinf(G) &= \inf\left\{ \entropysup_G(\ell) \mid \ell \in \calM^1(G)\right\} \\
\entropy_r &= \inf \left\{ \entropyinf(G) \mid \rk\pi_1 G = r \right\}
\end{align*}
\end{notation}

This this notation, Theorem~\ref{thm:main} can be restated as $\entropy_r > 0$.  As there are only finitely many finite graphs of any given rank without valence one or two vertices, the infimum in the definition of $\entropy_r$ is actually a minimum.

To simplify the analysis of the the matrix $A_{G,\ell}$ and the equilibrium measure $\mu$, it is necessary to have some assumptions about the structure of $G$.  In the previous section, we saw that having loop edges allowed for some connections to be made between the components of the vectors $\bu$ and $\bv$ that ultimately allowed for comparisons on components of $\mu$.  In this section we will show how to reduce proving Theorem~\ref{thm:main} to the case where every vertex has a loop edge.  Moreover, we will also show that we can assume the loops are short compared to the non-loop edges.

\begin{definition}\label{def:witness}
Suppose $G$ is a graph with $\rk \pi_1(G) \geq 3$.  We say a sequence of length functions $(\ell_i) \subset \calM^1(G)$ is \emph{witnessing} if $\entropysup_G(\ell_i) \to \entropyinf(G)$.  
\end{definition}

In other words, the sequence $(\ell_i)$ is witnessing if the largest possible entropy the restriction of $\ell_i$ to a proper subgraph of $G$ limits to the infimum over all proper subgraphs of rank $r$ graphs.  The next proposition shows how to modify a witnessing sequence in the setting where this infimum is equal to 0.
 
\begin{proposition}\label{prop:reduction}
Suppose $r \geq 3$ and $\entropy_r = 0$.  Then there is a graph $G$ with $\rk \pi_1(G) = r$, $\entropyinf(G) = 0$, and a witnessing sequence $(\ell_i) \subset \calM^1(G)$ such that:
\begin{enumerate}
\item\label{item:loop} every vertex of $G$ is incident to a loop edge,
\item\label{item:nonloop} there exists at least one non-loop edge, and
\item\label{item:short loops} for any non-loop edge $e$, there are loop edges $\gamma_1$ and $\gamma_2$ with $o(\gamma_1) = t(\gamma_1) = o(e)$ and $o(\gamma_2) = t(\gamma_2) = t(e)$, we have $\ell_i(\gamma_1), \ell_i(\gamma_2) \leq \frac{1}{4}\ell_i(e)$. 
\end{enumerate}
\end{proposition}

\begin{proof}
Fix $r \geq 3$ and suppose that $\entropy_r = 0$.  As there are only finitely many finite graphs with a fixed finite rank without valence 1 or 2 vertices, we may assume that we have a finite graph $G_1$ with $\rk\pi_1G_1 = r$, $\entropyinf(G_1) = 0$ and a witnessing sequence $(\ell_i) \subset \calM^1(G_1)$.  

First we will show how to arrange~\eqref{item:loop}.  Let $E_0$ be the set of edges of $G_1$ that are incident to a vertex of $G_1$ that is not incident to any loop edge.  If $E_0$ is empty, then we set $G' = G_1$.  Else, let $e_i$ be the edge in $E_0$ that minimizes $\ell_i(e_i)$.  As there are only finitely many edges, we may assume that the sequence $e_i$ is constant, i.e., $e_i = e$ for some edge $e \in E_0$.         
  
Let $G_2$ be the graph obtained by collapsing $e$.  There is an induced length function $\ell'_i$ in $G_2$ defined by $\ell'_i(e') = \ell_i(e')$ for each remaining edge.  Let $\alpha_i = \entropy_{G_2}(\ell'_i)$ so that $\hell_i = \alpha_i\ell'_i \in \calM^1(G_2)$. 

\begin{claim}\label{claim:alphi}
For any proper subgraph $H_2 \subseteq G_2$, there is a proper subgraph $H_1 \subseteq G_1$ such that \[ \entropy_{H_1}\left(\ell_i\big|_{H_1}\right) \leq \entropy_{H_2}\left(\ell'_i\big|_{H_2}\right) \leq 2\entropy_{H_1}\left(\ell_i\big|_{H_1}\right).\]
\end{claim}

\begin{proof}[Proof of Claim~\ref{claim:alphi}]
Let $H_2$ be a proper subgraph of $G_2$.  There is a unique proper subgraph $H_1 \subset G_1$ that contains $e$ and so that the collapse of $e$ results in $H_2$.  There is a bijection between cycles on $G_2$ and $G_1$ we denote by: $b \from \calC(G_2) \to \calC(G_1)$ which restricts to a bijection between cycles on $H_2$ and $H_1$.  

For any cycle $\gamma$ on $H_2$, by construction, we have that $\ell'_i(\gamma) \leq \ell_i(b(\gamma)) \leq 2\ell'_i(\gamma)$ as any cycle that crossed $e$ necessarily crossed a remaining edge that was longer.  Thus 
\begin{equation*}
\abs{\calC_{H_1,\ell_i}(t)} \leq \abs{\calC_{H_2,\ell'_i}(t)} \leq \abs{\calC_{H_1,\ell_i}(2t)}.
\end{equation*}  

The first of these inequalities immediately gives $\entropy_{H_1}\left(\ell_i\big|_{H_1}\right) \leq \entropy_{H_2}\left(\ell'_i\big|_{H_2}\right)$.  For the other we compute:
\begin{align*}
\entropy_{H_2}\left(\ell'_i\big|_{H_2}\right) &= \lim_{t \to \infty} \frac{1}{t} \log \abs{\calC_{H_2,\ell'_i}(t)}  \\
&\leq  \lim_{t \to \infty} \frac{1}{t} \log \abs{\calC_{H_1,\ell_i}(2t)} \\
&\leq  \lim_{\tau \to \infty} \frac{2}{\tau} \log \abs{\calC_{H_1,\ell_i}(\tau)} \\
& = 2\entropy_{H_1}\left(\ell_i\big|_{H_1}\right).\qedhere
\end{align*}
\end{proof}

In particular, we find that $\alpha_i = \entropy_{G_2}(\ell'_i) \geq \entropy_{G_1}(\ell_i) = 1$.  Combining this with the other inequality in Claim~\ref{claim:alphi} we find:
\begin{equation*}
\entropy_{H_2}\left(\hell_i|_{H_2}\right) = \frac{1}{\alpha_i}\entropy_{H_2}\left(\ell'_i|_{H_2}\right) \leq \entropy_{H_2}\left(\ell'_i|_{H_2}\right) \leq 2\entropy_{H_1}\left(\ell_i\big|_{H_1}\right).
\end{equation*}
Therefore we have that $\entropysup_{G_2}(\hell_i) \leq 2\entropysup_{G_1}(\ell_i)$ and thus $\entropysup_{G_2}(\hell_i) \to 0$.

We now repeat this procedure as needed replacing $G_1$ with $G_2$ so that the result is a graph $G'$ where every vertex is incident to a loop edge with $\entropyinf(G') = 0$, and a witnessing sequence $(\ell_i) \subset \calM^1(G')$.  This shows~\eqref{item:loop}

By Proposition~\ref{prop:rose}, $\entropyinf(\calR_r) > 0$, hence $G'$ must have a non-loop edge.  This shows~\eqref{item:nonloop}.

Lastly, to show how to arrange for~\eqref{item:short loops}, suppose there is a vertex where every incident loop edge has length more than $\frac{1}{4}$ of the length of every incident non-loop edge.  Then as we argued in showing~\eqref{item:loop}, we can contract the shortest of these non-loop edges and only change the length of cycles by a bounded amount.  An argument similar to the one above then gives the desired result.
\end{proof}

\section{Equilibrium Measure Bounds on Non-Loop Edges} \label{sec:lemma}

The main goal in this section is to prove an analogue of Lemma~\ref{lem:rose estimate} for non-loop edges that arise via the reduction procedure of the previous section.

\begin{lemma}\label{lem:estimate}
Let $G$ be a connected graph and fix a length function $\ell \in \calM^1(G)$ with equilibrium measure $\mu$.  Suppose that $e$ is a non-loop edge in $G$ that is incident at each of its vertices to a loop edge, denoted $\gamma_1$ and $\gamma_2$ respectively.  Then:   
\begin{equation*}
\exp(\ell(e))\mu(e) \leq 2\exp(\ell(\gamma_1) + \ell(\gamma_2))(\mu(\gamma_1) + \mu(\gamma_2))
\end{equation*}
\end{lemma}

\begin{proof}
Fix a length function $\ell \in \calM^1(G)$, let $A = A_{G,\ell}$ and fix vectors $\bu$ and $\bv$ such that $A\bv = \bv$, $\bu^T A = \bu^T$ and $\bu^T\bv = 1$.  Then we have that $\mu(e) = \bu(e) \bv(e)$ as explained in Section~\ref{subsec:entropy}.  Denote the edges at $o(e)$ other than $e$ and $\gamma_1$ by $e_1,\ldots,e_n$ and denote the edges at $t(e)$ other than $\bar{e}$ and $\gamma_2$ by $e'_1,\ldots,e'_m$.  See Figure~\ref{fig:estimate}.  Throughout the calculations below we will make use of the fact that $\bu(e') = \bu(\bar{e}')$ for any edge $e'$ and similarly for $\bv$ and $\mu$.

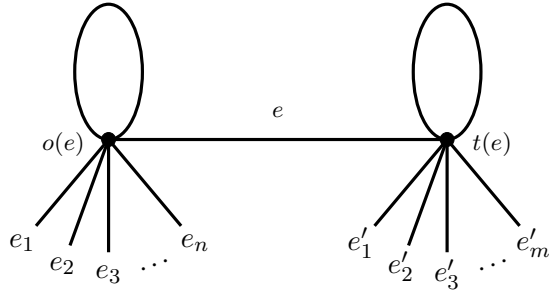
\begin{figure}[ht]
\centering
\begin{tikzpicture}[scale=1.5]
\draw[very thick] (-1.5,0) -- (1.5,0) node[pos=0,label={[xshift=-0.6cm,yshift=-0.4cm]\footnotesize $o(e)$},inner sep=0pt] (o) {} node[pos=0.5,label=above:{\footnotesize  $e$}] {} node[pos=1,label={[xshift=0.6cm,yshift=-0.4cm]\footnotesize $t(e)$},inner sep=0pt] (t) {};
\draw[very thick] (o) arc (-90:270: 0.3cm and 0.6cm);
\draw[very thick] (t) arc (-90:270: 0.3cm and 0.6cm);
\foreach \a in {-75,-70,-65} {
	\draw[white] (o) -- ++(\a:1cm) node[black,pos=1.2] {$\cdot$}; 
	\draw[white] (t) -- ++(\a:1cm) node[black,pos=1.2] {$\cdot$}; 
}
\foreach \a [count=\ai] in {-130,-110,-90,-50}{
	\ifnum\ai<4
		\draw[very thick] (o) -- ++(\a:1cm) node[pos=1.2] {$e_\ai$}; 
		\draw[very thick] (t) -- ++(\a:1cm) node[pos=1.2] {$e'_\ai$}; 
	\else
		\draw[very thick] (o) -- ++(\a:1cm) node[pos=1.2] {$e_n$}; 
		\draw[very thick] (t) -- ++(\a:1cm) node[pos=1.2] {$e'_m$}; 
	\fi
}
\filldraw (o) circle [radius=1.75pt];
\filldraw (t) circle [radius=1.75pt];
\end{tikzpicture}
\caption{The set up in Lemma~\ref{lem:estimate}.}\label{fig:estimate}
\end{figure}

As $\bu^T A = \bu^T$, we find that:
\begin{equation}\label{eq:u_e}
\bu(e) = 2\exp(-\ell(\gamma_1))\bu(\gamma_1) + \sum_{i=1}^n \exp(-\ell(e_i))\bu(e_i).
\end{equation}
We will denote the right-hand side of~\eqref{eq:u_e} by $X_\bu$ so that $\bu(e) = X_\bu$.  Similarly, as $A\bv = \bv$, we find that:
\begin{equation}\label{eq:v_e}
\bv(e) = \exp(-\ell(e))\left(2\bv(\gamma_2) + \sum_{i=1}^m \bv(e'_i)\right).
\end{equation}
We denote $2\bv(\gamma_2) + \sum_{i=1}^n \bv(e_i)$ by $X_\bv$ so that $\bv(e) = \exp(-\ell(e))X_\bv$.  Thus combining \eqref{eq:u_e} and \eqref{eq:v_e} we have
\begin{equation}\label{eq:mu_e}
\mu(e) = \bu(e)\bv(e) = \exp(-\ell(e))X_\bu X_\bv.
\end{equation}

We record the same calculations for $\bar{e}$ and find:
\begin{align*}
\bu(\bar{e}) & = 2\exp(-\ell(\gamma_2))\bu(\gamma_2) + \sum_{i=1}^m \exp(-\ell(e'_i))\bu(e'_i) = Y_\bu \\
\bv(\bar{e}) & = \exp(-\ell(e))\left(2\bv(\gamma_1) + \sum_{i=1}^n \bv(e_i)\right) = \exp(-\ell(e))Y_\bv
\end{align*}
Thus again we have that:
\begin{equation}\label{eq:mu_bare}
\mu(\bar{e}) = \exp(-\ell(e))Y_\bu Y_\bv.
\end{equation}

Next, we look at $\mu(\gamma_1)$ and $\mu(\gamma_2)$. We can compute them as above and compare these expressions to those for $\mu(e)$ and $\mu(\bar{e})$.  First we look at $\mu(\gamma_1) = \bu(\gamma_1)\bv(\gamma_1)$ and estimate these quantities.
\begin{align}
\bu(\gamma_1) &= \exp(-\ell(\gamma_1))\bu(\gamma_1) + \exp(-\ell(e))\bu(e) + \sum_{i=1}^n \exp(-\ell(e_i))\bu(e_i) \geq \frac{X_\bu}{2}\label{eq:u_gamma1} \\
\bv(\gamma_1) &= \exp(-\ell(\gamma_1))\left(\bv(\gamma_1) + \bv(e) + \sum_{i=1}^n \bv(e_i)\right) \geq \exp(-\ell(\gamma_1))\frac{Y_\bv}{2}\label{eq:v_gamma1}
\end{align}
In a similar manner, we also can estimate that $\bu(\gamma_2) \geq \frac{Y_\bu}{2}$ and $\bv(\gamma_2) \geq \exp(-\ell(\gamma_2))\frac{X_\bv}{2}$.  

Combining \eqref{eq:u_gamma1} and \eqref{eq:v_gamma1} together with these related inequalities, we find:
\begin{align*}
\mu(\gamma_1) + \mu(\gamma_2) & = \bu(\gamma_1)\bv(\gamma_1) + \bu(\gamma_2)\bv(\gamma_2) \\
& \geq \exp(-\ell(\gamma_1))\frac{X_\bu Y_\bv}{4} + \exp(-\ell(\gamma_2))\frac{X_\bv Y_\bu}{4} \\
& \geq \frac{1}{4}\exp(-\ell(\gamma_1)-\ell(\gamma_2))\left(X_\bu Y_\bv + X_\bv Y_\bu\right)
\end{align*}
As $\mu(e) = \mu(\bar{e})$, we have that $X_\bu X_\bv = Y_\bu Y_\bv$ by \eqref{eq:mu_e} and \eqref{eq:mu_bare}.  Hence $X_\bv = \frac{Y_\bu Y_\bv }{X_\bu}$ and $Y_\bv = \frac{X_\bu X_\bv}{Y_\bu}$.  Continuing the above inequalities, we find:
\begin{align*}
\mu(\gamma_1) + \mu(\gamma_2) &\geq \frac{1}{4}\exp(-\ell(\gamma_1)-\ell(\gamma_2))\left(X_\bu\left(\frac{X_\bu X_\bv}{Y_\bu}\right) + Y_\bu \left(\frac{Y_\bu Y_\bv }{X_\bu}\right)\right) \\
& = \frac{1}{4}\exp(-\ell(\gamma_1)-\ell(\gamma_2))\left(X_\bu X_\bv\left(\frac{X_\bu}{Y_\bu}\right) + Y_\bu Y_\bv \left(\frac{Y_\bu}{X_\bu}\right)\right) \\
& =  \frac{X_\bu X_\bv}{4}\exp(-\ell(\gamma_1)-\ell(\gamma_2))\left(\frac{X_\bu}{Y_\bu} + \frac{Y_\bu}{X_\bu}\right) \\
& \geq \frac{X_\bu X_\bv}{2} \exp(-\ell(\gamma_1) - \ell(\gamma_2)).
\end{align*}
The last inequality holds as $x + \frac{1}{x} \geq 2$ for $x \in (0,\infty)$.  Rewriting, we have:
\begin{equation}\label{eq:xx < gg}
X_\bu X_\bv \leq 2\exp(\ell(\gamma_1) + \ell(\gamma_2))(\mu(\gamma_1) + \mu(\gamma_2))
\end{equation}

The inequality \eqref{eq:xx < gg} together with \eqref{eq:mu_e} gives:
\begin{equation*}
\exp(\ell(e))\mu(e) = X_\bu X_\bv  \leq 2\exp(\ell(\gamma_1) + \ell(\gamma_2))(\mu(\gamma_1) + \mu(\gamma_2))
\end{equation*}
as desired.
\end{proof}

\section{The Proof of Theorem~\ref{thm:main}}\label{sec:proof}

In this section, we assemble the proof of Theorem~\ref{thm:main} and show that $\entropy_r > 0$.  (See Notation~\ref{notation:entropy}.) To this end, suppose the theorem is false, that is $\entropy_r = 0$ for some fixed $r \geq 3$.  Let $G$ be the graph with witnessing sequence $(\ell_i) \subset \calM^1(G)$ as given by Proposition~\ref{prop:reduction}.  Fix a non-loop edge $e$ in $G$ and loop edges, $\gamma_{1}, \gamma_{2}$, incident to the vertices $o(e)$ and $t(e)$ respectively. Recall that by~\eqref{item:short loops} in Proposition \ref{prop:reduction}, we can assume that 
\[ \ell_i (\gamma_j) \leq \frac{1}{4} \ell_i (e), \text{ for } j= 1,2. \]

For each $i$, let $\psi_t^{(i)}$ be the linear time blow-up of $\ell_i$ along $e$ and let $j_i \from \RR_{\geq 0} \to \RR_{\geq 0}$ be the corresponding scaling function.  Denote the equilibrium measure at $\psi_t^{(i)}$ by $\mu^{(i)}_t$.  Let $G' = G - e$.  By Proposition~\ref{prop:estimate} we have that:
\begin{equation}\label{eq:witnessing entropy}
\entropy_{G'}\left(\ell_i\big|_{G'}\right) = 1 - \int_0^\infty \abs{j'_i(t)} \, dt = 1 - \int_0^\infty \frac{\mu_t^{(i)}(e)}{\sum_{e' \neq e} \ell_i(e')\mu_t^{(i)}(e')} \, dt.
\end{equation}

By Lemma~\ref{lem:estimate} and the definition of $\psi^{(i)}_t$, we find that:
\begin{align*}
\int_0^\infty \frac{\mu_t^{(i)}(e)}{\sum_{e' \neq e} \ell(e')\mu_t^{(i)}(e')} \, dt & \leq \int_0^\infty \frac{2\exp(\psi^{(i)}_t(\gamma_1) + \psi^{(i)}_t(\gamma_2) - \psi^{(i)}_t(e))(\mu_t^{(i)}(\gamma_1) + \mu_t^{(i)}(\gamma_2))}{\sum_{e' \neq e} \ell_i(e')\mu_t^{(i)}(e')} \, dt \\ 
&=  \int_0^\infty \frac{2\exp(j_i(t)\ell_i(\gamma_1) + j_i(t)\ell_i(\gamma_2) - \ell_i(e) - t)(\mu_t^{(i)}(\gamma_1) + \mu_t^{(i)}(\gamma_2))}{\sum_{e' \neq e} \ell_i(e')\mu_t^{(i)}(e')} \, dt
\end{align*}

Let $m^{(i)} = \min\{ \ell_i(\gamma_1), \ell_i(\gamma_2) \}$.  Then the final integral in the above inequality can be rewritten, using the fact that 
\[ m^{(i)}(\mu_t^{(i)}(\gamma_1) + \mu_t^{(i)}(\gamma_2)) \leq \sum_{e' \neq e} \mu_t^{(i)}(e') \ell_i(e'), \]   as follows:
\begin{multline*}
\int_0^\infty \frac{2\exp(j_i(t)\ell_i(\gamma_1) + j_i(t)\ell_i(\gamma_2) - \ell_i(e) - t)(\mu_t^{(i)}(\gamma_1) + \mu_t^{(i)}(\gamma_2))}{\sum_{e' \neq e} \ell_i(e')\mu_t^{(i)}(e')} \, dt \\ = 2\int_0^\infty \frac{\exp(j_i(t)\ell_i(\gamma_1) + j_i(t)\ell_i(\gamma_2) - \ell_i(e) - t)}{m^{(i)}}\frac{m^{(i)}(\mu_t^{(i)}(\gamma_1) + \mu^{(i)}(\gamma_2))}{\sum_{e' \neq e} \ell_i(e')\mu_t^{(i)}(e')} \, dt \\
\leq 2\int_0^\infty \frac{\exp(j_i(t)\ell_i(\gamma_1) + j_i(t)\ell_i(\gamma_2) - \ell_i(e) - t)}{m^{(i)}} \, dt.
\end{multline*}

As $j_i(t)\ell_i(\gamma_1) \leq \ell_i(\gamma_1) \leq \frac{1}{4}\ell_i(e)$, we have that $\ell_i(e) - j_i(t)\ell_i(\gamma_1) - j_i(t)\ell(\gamma_2) \geq \frac{1}{2}\ell_i(e)$.  Therefore, we find:
\begin{equation*}
2\int_0^\infty \frac{\exp(j_i(t)\ell_i(\gamma_1) + j_i(t)\ell_i(\gamma_2) - \ell_i(e) - t)}{m^{(i)}} \, dt \leq \frac{2\exp(-\ell_i(e)/2)}{m^{(i)}}\int_0^\infty \exp(-t) \, dt = \frac{2\exp(-\ell_i(e)/2)}{m^{(i)}}
\end{equation*}  

This now gives:
\begin{equation*}
\entropy_{G'}\left(\ell_i\big|_{G'}\right) = 1 - \int_0^\infty \frac{\mu_t^{(i)}(e)}{\sum_{e' \neq e} \ell_i(e')\mu_t^{(i)}(e')} \, dt \geq 1 - \frac{2\exp(-\ell_i(e)/2)}{m^{(i)}}.
\end{equation*}
Since $\entropy_{G'}\left(\ell_i\big|_{G'}\right) \to \entropyinf(G) = 0$ as $i \to \infty$ by assumption, we must have that $\frac{2\exp(-\ell_i(e)/2)}{m^{(i)}} \to 1$ as $i \to \infty$.  In particular, for large enough $i$, we find that $m^{(i)} \leq 3\exp(-\ell_i(e)/2)$.  

Therefore, the edges $e, \gamma_1, \gamma_2$ form a sub-barbell $\calB \subset G$ so that one of $\gamma_1, \gamma_2$ has length at most $3 \exp(- \ell_i(e)/2)$ and so that the other self-loop has length at most $\frac{1}{4} \ell_i (e)$. Lemma~\ref{lem:barbell estimate} applies and we conclude that for all sufficiently large $i$ we have $\entropy_{\calB}\left(\ell_i \big|_{\calB}\right) \ge \frac{1}{5}$ and hence $\entropysup_{G} (\ell_i) \geq \frac{1}{5}$ as well.
This is a contradiction.

\bibliographystyle{acm}
\bibliography{subgraph}

@article {ar:ACR23,
    AUTHOR = {Aougab, Tarik and Clay, Matt and Rieck, Yo'av},
     TITLE = {Thermodynamic metrics on outer space},
   JOURNAL = {Ergodic Theory Dynam. Systems},
  FJOURNAL = {Ergodic Theory and Dynamical Systems},
    VOLUME = {43},
      YEAR = {2023},
    NUMBER = {3},
     PAGES = {729--793},
      ISSN = {0143-3857,1469-4417},
   MRCLASS = {20F65 (20E05 20E36 57K20)},
  MRNUMBER = {4544143},
       DOI = {10.1017/etds.2021.165},
       URL = {https://doi.org/10.1017/etds.2021.165}
}

@article {ar:Parlier24,
    AUTHOR = {Parlier, Hugo},
     TITLE = {A shorter note on shorter pants},
   JOURNAL = {Bull. Lond. Math. Soc.},
  FJOURNAL = {Bulletin of the London Mathematical Society},
    VOLUME = {56},
      YEAR = {2024},
    NUMBER = {4},
     PAGES = {1483--1487},
      ISSN = {0024-6093,1469-2120},
   MRCLASS = {57K20 (30F60 32G15)},
  MRNUMBER = {4743819},
MRREVIEWER = {Athanase\ Papadopoulos},
       DOI = {10.1112/blms.13007},
       URL = {https://doi.org/10.1112/blms.13007},
}

@incollection {ar:Bers74,
    AUTHOR = {Bers, Lipman},
     TITLE = {Spaces of degenerating {R}iemann surfaces},
 BOOKTITLE = {Discontinuous groups and {R}iemann surfaces ({P}roc. {C}onf.,
              {U}niv. {M}aryland, {C}ollege {P}ark, {M}d., 1973)},
    SERIES = {Ann. of Math. Stud.},
    VOLUME = {No. 79},
     PAGES = {43--55},
 PUBLISHER = {Princeton Univ. Press, Princeton, NJ},
      YEAR = {1974},
   MRCLASS = {30A46},
  MRNUMBER = {361051},
MRREVIEWER = {C.\ Earle},
}

@Article{ar:M96, 
AUTHOR = {Minsky, Yair}, 
TITLE = {Extremal length estimates and product regions in Teichm{\"u}ller space}, 
JOURNAL = {Duke Math. Journal}, 
VOLUME = {83}, 
PAGES = {249--286}, 
YEAR = {1996} 

}

@Article{ar:MM00, 
AUTHOR = {Masur, Howard and Minsky, Yair}, 
TITLE= {Geometry of the complex of curves, II: Hierachical structure}, 
JOURNAL = {Geom. Func. Anal.}, 
VOLUME = {10}, 
PAGES = {902--974}, 
YEAR = {2000}


}

@phdthesis{th:Buser80,
  title={Riemannsche Fl{\"a}chen und L{\"a}ngenspektrum vom trigonometrischen Standpunkt aus},
  author={Buser, Peter},
  year={1980},
  school={Mathematisches Institut der Univ.}
}

@article {ar:schmutz94,
    AUTHOR = {Schmutz, Paul},
     TITLE = {Systoles on {R}iemann surfaces},
   JOURNAL = {Manuscripta Math.},
  FJOURNAL = {Manuscripta Mathematica},
    VOLUME = {85},
      YEAR = {1994},
    NUMBER = {3-4},
     PAGES = {429--447},
      ISSN = {0025-2611,1432-1785},
   MRCLASS = {57M50 (30F60 53C23)},
  MRNUMBER = {1305753},
MRREVIEWER = {Athanase\ Papadopoulos},
       DOI = {10.1007/BF02568209},
       URL = {https://doi.org/10.1007/BF02568209},
}

@article{ar:BS92,
  title={Symmetric pants decompositions of Riemann surfaces},
  author={Buser, Peter and Sepp{\"a}l{\"a}, Mika},
  journal = {Duke Math. Journal}, 
  volume = {67}, 
  number = {1}, 
  pages = {39--55},
  year={1992}
}

@article {ar:Brock03,
    AUTHOR = {Brock, Jeffrey F.},
     TITLE = {The {W}eil-{P}etersson metric and volumes of 3-dimensional
              hyperbolic convex cores},
   JOURNAL = {J. Amer. Math. Soc.},
  FJOURNAL = {Journal of the American Mathematical Society},
    VOLUME = {16},
      YEAR = {2003},
    NUMBER = {3},
     PAGES = {495--535},
      ISSN = {0894-0347,1088-6834},
   MRCLASS = {32G15 (30F40 30F60 37F30)},
  MRNUMBER = {1969203},
MRREVIEWER = {Edward\ C.\ Taylor},
       DOI = {10.1090/S0894-0347-03-00424-7},
       URL = {https://doi.org/10.1090/S0894-0347-03-00424-7},
}

@article {ar:Bestvina,
    AUTHOR = {Bestvina, Mladen},
     TITLE = {A {B}ers-like proof of the existence of train tracks for free
              group automorphisms},
   JOURNAL = {Fund. Math.},
  FJOURNAL = {Fundamenta Mathematicae},
    VOLUME = {214},
      YEAR = {2011},
    NUMBER = {1},
     PAGES = {1--12},
      ISSN = {0016-2736,1730-6329},
   MRCLASS = {20E05 (20E36)},
  MRNUMBER = {2845630},
MRREVIEWER = {Katalin\ A.\ Bencsath},
       DOI = {10.4064/fm214-1-1},
       URL = {https://doi.org/10.4064/fm214-1-1},
}

@article{ar:BestvinaFeighn, 
title = {Hyperbolicity of the complex of free factors}, 
author = {Bestvina, Mladen and Feighn, Mark}, 
journal = {Adv. Math.}, 
volume = {256}, 
pages = {104--155}, 
year = {2014}

}

@article{ar:CullerVogtmann, 
title ={Moduli of graphs and automorphisms of free groups}, 
author = {Culler, Mark and Vogtmann, Karen}, 
journal = {Invent. Math.}, 
volume = {84}, 
number = {1}, 
pages= {91--119}, 
year = {1986} 
}

@article{ar:Kao, 
author = {Kao, Lien Yung}, 
title = {Pressure type metrics on spaces of metric graphs}, 
journal = {Geom. Dedicata},
volume = {187}, 
pages = {151--177}, 
year = {2017} 

}

@Article{ar:PS14,
  author        = {Pollicott, Mark and Sharp, Richard},
  journal       = {Geom. Dedicata},
  title         = {A {W}eil-{P}etersson type metric on spaces of metric graphs},
  year          = {2014},
  issn          = {0046-5755},
  pages         = {229--244},
  volume        = {172},
  doi           = {10.1007/s10711-013-9918-2},
  fjournal      = {Geometriae Dedicata},
  mrclass       = {05C25 (32G15 37C30)},
  mrnumber      = {3253781},
  mrreviewer    = {Muhammad Salman},
  url           = {http://dx.doi.org/10.1007/s10711-013-9918-2}
}

@Article{ar:PP90,
  author        = {Parry, William and Pollicott, Mark},
  journal       = {Ast{\'e}risque},
  title         = {Zeta functions and the periodic orbit structure of hyperbolic dynamics},
  year          = {1990},
  issn          = {0303-1179},
  number        = {187-188},
  pages         = {268},
  fjournal      = {Ast{\'e}risque},
  mrclass       = {58F20 (58F11 58F15)},
  mrnumber      = {1085356},
  mrreviewer    = {Nicola{\u\i} T. A. Haydn}
}

@Article{ar:McMullen15,
  author        = {McMullen, Curtis T.},
  journal       = {J. Topol.},
  title         = {Entropy and the clique polynomial},
  year          = {2015},
  issn          = {1753-8416},
  number        = {1},
  pages         = {184--212},
  volume        = {8},
  fjournal      = {Journal of Topology},
  mrclass       = {37B40 (05C20 05C31 05C50)},
  mrnumber      = {3335252},
  mrreviewer    = {Steven M. Pederson},
  url           = {https://doi.org/10.1112/jtopol/jtu022}
}

@Article{Kim2020,
  author   = {Wooyeon Kim and Seonhee Lim},
  journal  = {Discrete and Continuous Dynamical Systems},
  title    = {Notes on the values of volume entropy of graphs},
  year     = {2020},
  issn     = {1078-0947},
  number   = {9},
  pages    = {5117-5129},
  volume   = {40},
  doi      = {10.3934/dcds.2020221},
  url      = {https://www.aimsciences.org/article/id/49387793-1b74-4973-87a7-a1d8b1499295}
}

@Book{bk:Walters82,
  author     = {Walters, Peter},
  publisher  = {Springer-Verlag, New York-Berlin},
  title      = {An introduction to ergodic theory},
  year       = {1982},
  isbn       = {0-387-90599-5},
  series     = {Graduate Texts in Mathematics},
  volume     = {79},
  mrclass    = {28Dxx (54H20 58F11)},
  mrnumber   = {648108},
  mrreviewer = {M.\ A.\ Akcoglu},
  pages      = {ix+250}
  
}

@article{ar:Wol, 
title = {Thurston's {R}iemannian metric for {T}eichm{\"u}ller space}, 
author = {Wolpert, Scott},
journal = {Journal of Differential Geometry}, 
volume = {23}, 
number = {2}, 
pages = {143--174}, 
year = {1986} 
}

@article {ar:Mc,
    AUTHOR = {McMullen, Curtis T.},
     TITLE = {Thermodynamics, dimension and the {W}eil-{P}etersson metric},
   JOURNAL = {Invent. Math.},
  FJOURNAL = {Inventiones Mathematicae},
    VOLUME = {173},
      YEAR = {2008},
    NUMBER = {2},
     PAGES = {365--425},
      ISSN = {0020-9910,1432-1297},
   MRCLASS = {37F30 (30F35 30F60 32G15 37D35)},
  MRNUMBER = {2415311},
MRREVIEWER = {Zheng\ Huang},
       DOI = {10.1007/s00222-008-0121-2},
       URL = {https://doi.org/10.1007/s00222-008-0121-2},
}

@Book{Farb2012,
  author        = {Farb, Benson and Margalit, Dan},
  publisher     = {Princeton University Press, Princeton, NJ},
  title         = {A primer on mapping class groups},
  year          = {2012},
  isbn          = {978-0-691-14794-9},
  series        = {Princeton Mathematical Series},
  volume        = {49},
  mrclass       = {57M50 (20F36 20F65 57M07 57N05)},
  mrnumber      = {2850125},
  mrreviewer    = {Stephen P. Humphries},
  pages         = {xiv+472}
}

\end{document}